\documentclass[12pt]{amsart}

\usepackage{a4wide}
\usepackage{graphicx,eepic}
\usepackage{amsmath,amsfonts,amssymb}
\usepackage{times}
\usepackage{bm}

\newcommand{\la}{\lambda}
\newcommand{\be}{\beta}
\newcommand{\de}{\delta}
\newcommand{\De}{\Delta}

\newcommand{\ga}{\gamma}
\newcommand{\e}{\varepsilon}

\newcommand{\om}{\omega}

\newcommand{\Si}{\Sigma}
\newcommand{\BR}{\mathbb{R}}
\newcommand{\BQ}{\mathbb{Q}}
\newcommand{\BZ}{\mathbb{Z}}
\newcommand{\BN}{\mathbb{N}}

\newcommand{\U}{\mathcal{U}}
\newcommand{\m}{\mathfrak{m}}

\newcommand{\kl}{q_{\mathrm{KL}}}

\newcommand{\B}{\mathcal{B}}
\newcommand{\V}{\mathcal{V}}

\newtheorem{lemma}{Lemma}[section]

\newtheorem{prop}[lemma]{Proposition}
\newtheorem{thm}[lemma]{Theorem}
\newtheorem{cor}[lemma]{Corollary}
\theoremstyle{definition}

\theoremstyle{remark}
\newtheorem{rmk}[lemma]{Remark}

\numberwithin{equation}{section} \numberwithin{table}{section}

\title[Expansions in non-integer bases]
{Expansions in non-integer bases: \\
lower, middle and top orders}
\author{Nikita Sidorov}
\address{School of Mathematics, The University of Manchester, Oxford Road, Manchester M13 9PL, United Kingdom. E-mail:
sidorov@manchester.ac.uk}
\date{August 27, 2008}
\dedicatory{To the memory of Bill Parry} 

\begin{document}

\begin{abstract}
Let $q\in(1,2)$; it is known that each $x\in[0,1/(q-1)]$ has an expansion of the form $x=\sum_{n=1}^\infty a_nq^{-n}$ with
$a_n\in\{0,1\}$. It was shown in \cite{EJK} that if
$q<(\sqrt5+1)/2$, then each $x\in(0,1/(q-1))$ has a continuum of such expansions; however, if $q>(\sqrt5+1)/2$, then there exist infinitely many $x$ having a unique expansion \cite{GS}.

In the present paper we begin the study of parameters $q$ for which there exists $x$ having a fixed finite number $m>1$ of expansions in base~$q$. In particular, we show that if $q<q_2=1.71\dots$, then each $x$ has either 1 or infinitely many expansions, i.e., there are no such $q$ in $((\sqrt5+1)/2,q_2)$.

On the other hand, for each $m>1$ there exists $\ga_m>0$ such that for any $q\in(2-\ga_m,2)$, there exists $x$ which has exactly $m$ expansions in base~$q$.
\end{abstract}

\maketitle

\section{Introduction and summary}
\label{sec:intro}

Expansions of reals in non-integer bases have been studied since the late 1950s, namely, since the pioneering works by R\'enyi \cite{Re} and Parry \cite{Pa}. The model is as follows: fix $q\in(1,2)$ and call any 0-1 sequence $(a_n)_{n\ge1}$ an {\it expansion in base~$q$} for some $x\ge0$ if
\begin{equation}\label{beta}
x=\sum_{n=1}^\infty a_nq^{-n}.
\end{equation}
Note that $x$ must belong to $I_q:=[0,1/(q-1)]$ and that for each $x\in I_q$ there is always at least one way of obtaining the $a_n$, namely, via the greedy algorithm (``choose 1 whenever possible'') -- which until recently has been considered virtually the only option.

In 1990 Erd\H os {\it et al.} \cite{EJK} showed (among other things) that if $q<G:=(\sqrt5+1)/2\approx1.61803$, then each
$x\in(0,1/(q-1))$ has in fact $2^{\aleph_0}$ expansions of the
form~(\ref{beta}). If $q=G$, then each $x\in I_q$ has $2^{\aleph_0}$ expansions, apart from $x=nG\pmod 1$ for $n\in\BZ$, each of which has $\aleph_0$ expansions in base~$q$ (see \cite{SV} for a detailed study of the space of expansions for this case). However, if $q>G$, then although a.e. $x\in I_q$ has $2^{\aleph_0}$ expansions in base~$q$ \cite{amm}, there always exist (at least countably many) reals having a unique expansion -- see \cite{GS}.

Let $\U_q$ denote the set of $x\in I_q$ which have a unique
expansion in base~$q$. The structure of the set $\U_q$ is reasonably well understood; its main property is that $\U_q$ is countable if $q$ is ``not too far" from the golden ratio, and uncountable of Hausdorff dimension strictly between 0 and 1 otherwise. More precisely, let $\kl$ denote the
\textit{Komornik-Loreti constant} introduced in \cite{KL}, which is defined as the unique solution of the equation
$$
\sum_{1}^{\infty}\mathfrak{m}_{n}x^{-n}=1,
$$
where $\mathfrak{m}=(\mathfrak{m}_n)_0^\infty$ is the Thue-Morse sequence $\mathfrak{m}=0110\,\,1001\,\,1001\,\,0110\,\dots$, i.e., a fixed point of the morphism $0\to01,\ 1\to10$. The Komornik-Loreti constant is known to be the smallest $q$ for which $x=1$ has a unique expansions in base~$q$ (see \cite{KL}), and its numerical value is approximately  $1.78723$.\footnote{For the list of all constants used in the present paper, see Table~\ref{table} before the bibliography.}

It has been shown by Glendinning and the author in \cite{GS} that
\begin{enumerate}
\item $\U_q$ is countable if $q\in(G,\kl)$, and each unique expansion is eventually periodic;
\item $\U_q$ is a continuum of positive Hausdorff dimension if $q>\kl$.
\end{enumerate}

Let now $m\in\BN\cup\{\aleph_0\}$ and put
\begin{align*}
\B_m=\{&q\in(G,2): \exists x\in I_q\ \mbox{which has exactly}\ m\\
&\mbox{expansions in base $q$ of the form}~(\ref{beta})\}.
\end{align*}
It follows from the quoted theorem from \cite{GS} that $\B_1=(G,2)$, but very little has been known about $\B_m$ for $m\ge2$. The purpose of this paper is to begin a systematic study of these sets.

\begin{rmk}
It is worth noting that in \cite{E3} it has been shown that for each $m\in\BN$ there exists an uncountable set $E_m$ of $q$ such that the number $x=1$ has $m+1$ expansions in base~$q$. The set $E_m\subset(2-\e_m,2)$, where $\e_m$ is small. A similar result holds for $m=\aleph_0$.

Note also that a rather general way to construct numbers
$q\in(1.9,2)$ such that $x=1$ has two expansions in base~$q$, has been suggested in \cite{KL2}.
\end{rmk}

\section{Lower order: $q$ close to the golden ratio}

We will write $x\sim(a_1,a_2,\dots)_q$ if $(a_n)_{n\ge1}$ is an expansion of $x$ in base~$q$ of the form~(\ref{beta}).

\begin{thm}\label{trans} For any transcendental $q\in(G,\kl)$
we have the following dichotomy: each $x\in I_q$ has either a
unique expansion or a continuum of expansions in base~$q$.
\end{thm}
\begin{proof} We are going to exploit the idea of {\it branching} introduced in \cite{Sid}. Let $x\in I_q$ have at least two expansions of the form~(\ref{beta}); then there exists the smallest $n\ge0$ such
that $x\sim(a_1,\dots,a_n,a_{n+1},\dots)_q$ and
$x\sim(a_1,\dots,a_n,b_{n+1},\dots)_q$ with $a_{n+1}\neq b_{n+1}$.
We may depict this {\em bifurcation} as shown in
Fig.~\ref{fig:branching}.

\begin{figure}[t]
    \centering 
\centerline{
\begin{picture}(265,210)
\thicklines  \put(0,100){\line(1,0){100}}
\put(100,100){\line(0,-1){50}}
\put(100,100){\line(0,1){50}}\put(100,150){\line(1,0){120}}
\put(100,50){\line(1,0){100}}
\put(220,150){\line(0,1){30}} \put(220,150){\line(0,-1){30}}
\put(200,50){\line(0,1){30}} \put(200,50){\line(0,-1){30}}
\put(200,20){\line(1,0){50}}\put(200,80){\line(1,0){50}}
\put(220,120){\line(1,0){50}}\put(220,180){\line(1,0){50}}
\put(274,179){$\dots$}\put(274,119){$\dots$}
\put(254,79){$\dots$}\put(254,19){$\dots$}
\put(0,110){$a_1\,\,\,\,\, a_2\,\,\dots\,\, a_{n-1}\,\,\, a_n$}
\put(103,155){$a_{n+1}\,\,\,a_{n+2}\,\,\,\dots\,\,\dots\,\,
a_{n_2}$}
\put(103,55){$b_{n+1}\,\,\,b_{n+2}\,\,\,\dots\,\,\,b_{n_2'}$}
\put(-10,97){$x$}
\end{picture}}
\caption{Branching and bifurcations}
    \label{fig:branching}
\end{figure}
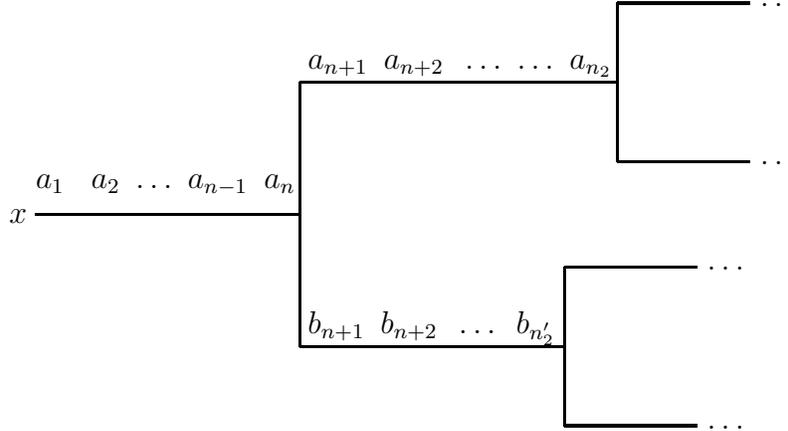

If $(a_{n+1},a_{n+2},\dots)_q$ is not a unique expansion, then there exists $n_2>n$ with the same property, etc. As a result, we obtain a subtree of the binary tree which corresponds to the set of all expansions of $x$ in base~$q$, which we call the {\it branching tree of $x$}. It has been shown in \cite[Theorem~3.6]{Sid} that if $q\in(G,\kl)$ that for for all $x$, except, possibly, a countable set, the branching tree is in fact the full binary tree and hence $x$ has $2^{\aleph_0}$ expansions in base~$q$; the issue is thus about these exceptional $x$'s.

Note that for $x$ to have at most countably many expansions in
base~$q$, its branching tree must have at least two branches which do not bifurcate. In other words, there exist two expansions of $x$ in base~$q$, $(a_n)_{n\ge1}$ and $(b_n)_{n\ge1}$ such that $(a_k,a_{k+1},\dots)$ is a unique expansion and so is $(b_j,b_{j+1},\dots)$ for some $k,j\in\BN$.

Without loss of generality, we may assume $j=k$, because the shift of a unique expansion is known to be a unique expansion \cite{GS}. Hence
\begin{align*}
x&=\sum_{i=1}^{k} a_iq^{-i} + q^{-k}r_k(q)\\
 &=\sum_{i=1}^{k} b_iq^{-i} + q^{-k}r'_k(q),
\end{align*}
where $r_k(q),r_k'(q)\in\U_q$ and $r_k(q)\neq r_k'(q)$. (If they are equal, then $q$ is obviously algebraic.) Since each unique expansion for $q\in(G, \kl)$ is eventually periodic (\cite[Proposition~13]{GS}), we have $\U_q\subset\BQ(q)$, whence the equation \begin{equation}\label{eq:aux1}
\sum_{i=1}^k (a_i-b_i)q^{-i}=q^{-k}(r_k(q)-r_k'(q))
\end{equation}
implies that $q$ is algebraic, unless (\ref{eq:aux1}) is an
identity. Assume it is an identity for some $q$; then it is an
identity for all $q>1$, because $r_k(q)=\pi(q)+\rho(q)/(1-q^{-r})$ and $r'_k(q)=\pi'(q)+\rho'(q)/(1-q^{-r'})$, where
$\pi,\pi',\rho,\rho'$ are polynomials.

Let $j=\min\,\{i\ge1:a_i\neq b_i\}<k$. We multiply (\ref{eq:aux1}) by $q^j$ and get
\[
a_j-b_j+\sum_{i=j+1}^k(a_i-b_i)q^{j-i}\equiv
q^{j-k}(r_k(q)-r'_k(q)),
\]
which is impossible, since $q\to+\infty$ implies $a_j-b_j=0$, a contradiction.
\end{proof}

The next question we are going to address in this section is
finding the smallest element of $\B_2$. Let $q\in\B_m$ and denote by $\U_q^{(m)}$ the set of $x\in I_q$ which have $m$ expansions in base~$q$. Firstly, we give a simple characterization of the set $\B_2$:

\begin{lemma}\label{crit}
A number $q\in(G,2)$ belongs to $\B_2$ if and only if
$1\in\U_q-\U_q$.
\end{lemma}
\begin{proof}1. Let $q\in\B_2$; then there exists $x$ having exactly two expansions in base~$q$. Without loss of generality we may assume that there exist two expansions of $x$, with $a_1=0$ and with $a_1=1$. (Otherwise we shift the expansion of $x$ until we obtain $x'$ having this property.) Note that
$x\in\bigl[\frac1q,\frac1{q(q-1)}\bigr]=:J_q$ -- the interval which is called the switch region in \cite{DK}.

Conversely, if $x\in J_q$, then it has a branching at
$n=1$. Since $x$ has only two different expansions in base~$q$, both shifts of $x$, namely, $qx$ (for $a_1=0$) and $qx-1$ (for $a_1=1$), must belong to $\U_q$, whence $1\in\U_q-\U_q$.

\smallskip\noindent 2. Let $y\in\U_q$ and $y+1\in\U_q$.
We claim that $x:=(y+1)/q$ belongs to $\U_q^{(2)}$.
Note that $y\in\U_q$ implies $y\not\in J_q$, whence $y<1/q$, because if $y$ were greater than $1/(q(q-1))$, we would have
$y+1>\frac{q^2-q+1}{q-1}>\frac1{q-1}$.

Thus, $y<1/q$, whence $x\in J_q$, because $y+1<1/(q-1)$.  Since $x\in J_q$, it has at least two different expansions in base~$q$, with $a_1=0$ and $a_1=1$, and shifting each of them yields $qx=y+1$ and $qx-1=y$, both having unique expansions. Hence there are only two possible expansions of $x$, i.e., $x\in\U_q^{(2)}$.
\end{proof}

This criterion, simple as it is, indicates the difficulties one faces when dealing with $\B_2$ as opposed to the unique expansions; at first glance, it may seem rather straightforward to verify whether if a number $x$ has a unique expansion, then so does $x+1$ -- but this is not the case.

The reason why this is actually hard is the fact that ``typically'' adding 1 to a number alters the tail of its greedy expansion (which, of course, coincides with its unique expansion if $x\in\U_q$) in a completely unpredictable manner -- so there is no way of telling whether $x+1$ belongs to $\U_q$ as well.

Fortunately, if $q$ is sufficiently small, the set of unique
expansions is very simple, and if $q$ is close to 2, then $\U_q$ is large enough to satisfy $\U_q-\U_q=[-1/(q-1),1/(q-1)]$ -- see Section~\ref{sec:top}.

\begin{lemma}\label{lem:gs}
Let $G<q\le q_f$; then any unique expansion belongs to the set $\{0^k(10)^\infty, 1^k(01)^\infty,\linebreak 0^\infty,1^\infty\}$ with $k\ge0$.
\end{lemma}
\begin{proof}If $x\in \De_q:=((2-q)/(q-1),1)$, then, by \cite[Section~4]{GS}, each unique expansion for this range of $q$ is either $(10)^\infty$ or $(01)^\infty$. If $x\in I_q\setminus\De_q$, then any unique expansion is of the form $1^k\e$ or $0^k\e$, where $\e$ is a unique expansion of some $y\in\De_q$ (\cite[Corollary~15]{GS}).
\end{proof}

\begin{prop}\label{thm1} The smallest element of $\B_2$ is $q_2$, the appropriate root of
\begin{equation}\label{q2}
x^4=2x^2+x+1,
\end{equation}
with a numerical value $1.71064$. Furthermore,
$\B_2\cap(G,q_f)=\{q_2\}$.
\end{prop}
\begin{proof} Let $q_f$ be the cubic unit which satisfies
\begin{equation}\label{qf}
x^3=2x^2-x+1,\quad q_f\approx1.75488\dots
\end{equation}
We first show that $q_f\in\B_2$. By Lemma~\ref{crit}, it suffices to produce $y\in\U_q$ such that $y+1\in\U_q$ as well. Note that $q_f$ satisfies $x^4=x^3+x^2+1$ (together with $-1$); put $y\sim(0000010101\dots)_{q_f}$. Then
$y+1\sim(11010101\dots)_{q_f}$, both unique expansions by
Lemma~\ref{lem:gs}.

Hence $\inf\B_2\le q_f$. This makes our search easier, because by Lemma~\ref{lem:gs}, each unique expansion for $q\in(G,q_f)$
belongs to the set $\{0^k(10)^\infty, 1^k(01)^\infty,
0^\infty,1^\infty\}$ with $k\ge0$.

Let us show first that the two latter cases are impossible for
$q\in(G,q_f)$. Indeed, if $x\sim(10^\infty)_q$ had exactly two
expansions in base~$q$, then the other expansion would be of the
form $(01^k(01)^\infty)_q$, which would imply
$1=1/q+1/q^2+\dots+1/q^k+1/q^{k+2}+1/q^{k+4}+\dots$ with $k\ge1$. If
$k\ge2$, then $1<1/q+1/q^2+1/q^4$, i.e., $q>q_f$; $k=1$ implies
$q=G$. The case of the tail $1^\infty$ is completely
analogous

To simplify our notation, put $\la=q^{-1}\in(1/q_f,1/G)$. So let
$x\sim (0^{\ell-1}(10)^\infty)_q$ and
$x+1\sim(1^{k-1}(01)^\infty)_q$, both in $\U_q$, with $\ell\ge1,
k\ge1$. Then we have
\[
1+\frac{\la^{\ell}}{1-\la^2}=\frac{\la-\la^{k-1}}{1-\la}+
\frac{\la^{k-1}}{1-\la^2}.
\]
Simplifying this equation yields
\begin{equation}\label{main}
\la^{\ell}+\la^k=2\la^2+\la-1.
\end{equation}
In view of symmetry, we may assume $k\ge\ell$.

\medskip\noindent\textbf{Case 1:} $\ell=1$. This implies
$\la^k=2\la^2-1$, whence $2\la^2-1>0$, i.e., $\la>1/\sqrt2>1/G$. Thus, there are no solutions of (\ref{main}) lying in $(1/q_f,1/G)$
for this case.

\medskip\noindent\textbf{Case 2:} $\ell=2$. Here
$\la^k=\la^2+\la-1>0$, whence $\la>1/G$. Thus, there are no
solutions here either.

\medskip\noindent\textbf{Case 3:} $\ell=3$. We have
\begin{equation}\label{ell3}
\la^k = -\la^3+2\la^2+\la-1.
\end{equation}
Note that the root of (\ref{ell3}) as a function of $k$ is
decreasing. For $k=3$ the root is above $1/G$, for $k=4$ it is
exactly $1/G$. For $k=5$ the root of (\ref{ell3}) satisfies
$x^5=-x^3+2x^2+x-1$, which can be factorized into $x^4+x^3+2x^2=1$,
i.e., the root is exactly $1/q_2$.

Finally, for $k=6$ the root satisfies $x^6=-x^3+2x^2+x-1$, which factorizes into $x^3-x^2+2x-1=0$, i.e., $\la=1/q_f$. For $k>6$ the root of (\ref{ell3}) lies outside the required range.

\medskip\noindent\textbf{Case 4:} $\ell=4, k\in\{4,5\}$. For $k=4$ the root of $2x^4=2x^2+x-1$ is $0.565\ldots<1/q_f=0.569\ldots$. If $k=5$, then the root is $0.543\dots$, i.e., even smaller. Hence there are no appropriate solutions of (\ref{main}) here.

\medskip\noindent\textbf{Case 5:}
If $\ell\ge5$ and $k\ge5$, then the LHS of (\ref{main}) is less than $2G^{-5}<0.2$, whereas the RHS is greater than
$2q_f^{-2}+q_f^{-1}-1>0.21$, whence there are no solutions of
(\ref{main}) in this case. If $\ell=4,k\ge6$, then, similarly,
$\la^k+\la^\ell\le\la^4+\la^6<G^{-4}+G^{-6}<0.202$.

\bigskip

Thus, the only case which produces a root in the required range is Case~3, which yields $1/q_2$. Hence
\begin{equation}\label{eq:capb2}
(G,q_f)\cap\B_2=\{q_2\}.
\end{equation}
\end{proof}

\begin{rmk} Let $q=q_2$ and let $y\sim(0000(10)^\infty)_{q_2}
\in\U_q$ and $y+1\sim (11(01)^\infty)_{q_2}\in\U_q$. We thus see that in this case the {\em tail} of the expansion does change, from $(10)^\infty$ to $(01)^\infty$. (Not the {\em period}, though!) Also, the proof of Lemma~\ref{crit} allows us to construct $x\in\U_{q_2}^{(2)}$ explicitly, namely,
$x\sim(011(01)^\infty)_{q_2}\sim(10000(10)^\infty)_{q_2}$, i.e., $x\approx0.64520$.
\end{rmk}

A slightly more detailed study of equation~(\ref{main}) shows that it has only a finite number of solutions $\la\in(1/\kl,1/G)$. In order to construct an infinite number of $q\in\B_2\cap(q_f,\kl)$, one thus needs to consider unique expansions with tails different from $(01)^\infty$:

\begin{prop}The set $\B_2\cap(q_f,\kl)$ is infinite countable.
\end{prop}
\begin{proof}We are going to develop the idea we used to show that $q_f\in\B_2$. Namely, let $(q_f^{(n)})_{n\ge1}$ be the sequence of algebraic numbers specified by their greedy expansions of 1:
\begin{align*}
q_f^{(1)}&:1\sim (11\,\,0^\infty)_{q_f^{(1)}}=G,\\
q_f^{(2)}&: 1\sim (1101\,\,0^\infty)_{q_f^{(2)}}=q_f\\
q_f^{(3)}&: 1\sim (1101\,\, 0011\,\,0^\infty)_{q_f^{(3)}}\\
&\vdots\\
q_f^{(n)}&: 1\sim (\m_1,\dots,\m_{2^n}\,\,0^\infty)_{q_f^{(n)}},
\end{align*}
where $(\m_n)$ is the Thue-Morse sequence -- see Introduction. It is obvious that $q_f^{(n)}\nearrow\kl$. We now define the sequence $z_n$ as follows:
\[
z_n\sim (0^{2^n}(\m_{2^{n-1}+1}\dots\m_{2^n})^\infty)_{q_f^{(n)}},
\]
whence
\[
z_n+1\sim (\m_1,\dots,\m_{2^{n-1}}
(\m_{2^{n-1}+1}\dots\m_{2^n})^\infty)_{q_f^{(n)}}.
\]
\cite[Proposition~9]{GS} implies that $z_n\in\U_{q_f^{(n)}}$ and $z_n+1\in\U_{q_f^{(n)}}$, whence by Lemma~\ref{crit}, $q_f^{(n)}\in\B_2$ for all $n\ge2$.
\end{proof}

\begin{lemma}\label{r1} We have $\B_m\subset\B_2$ for any natural $m\ge3$.
\end{lemma}
\begin{proof}
If $q\in\B_m$ for some natural $m\ge3$, then the branching argument immediately implies that there exists $x\in\U_q^{(m')}$, with $1<m'<m$. Hence, by induction, there exists $x'\in\U_q^{(2)}$. Therefore, $\B_m\subset\B_2$ for all $m\in\BN\setminus\{1\}$.
\end{proof}

Our next result shows that a weaker analogue of Theorem~\ref{trans} holds without assuming $q$ being transcendental, provided $q<q_f$.

\begin{thm}\label{thm2}
For any $q\in(G,q_2)\cup(q_2,q_f)$, each $x\in I_q$ has either a unique expansion or infinitely many expansions of the
form~(\ref{beta}) in base~$q$. Here $G=\frac{1+\sqrt5}2$ and $q_2$ and $q_f$ are given by (\ref{q2}) and (\ref{qf}) respectively.
\end{thm}
\begin{proof} It follows immediately from Proposition~\ref{thm1}, Lemma~\ref{r1} and relation~(\ref{eq:capb2}) that $\B_m\cap(G,q_f)\subset\{q_2\}$ for all $m\in\BN\setminus\{1\}$.
\end{proof}

\begin{cor}For $q\in(G,q_2)\cup(q_2,q_f)$ each $x\in J_q$ has
infinitely many expansions in base~$q$.
\end{cor}
\begin{proof}It suffices to recall that each $x\in J_q$ has at least two expansions in base~$q$ and apply Theorem~\ref{thm2}.
\end{proof}

It is natural to ask whether the claim of Theorem~\ref{thm2} can be strengthened in the direction of getting rid of $q\in\B_{\aleph_0}$ so we could claim that a stronger version of Theorem~\ref{trans} holds for $q<q_f$. It turns out that the answer to this question is negative.

Notice first that $\B_{\aleph_0}\not\subset\B_2$, since $G\in\B_{\aleph_0}\setminus\B_2$. Our goal is to
show that in fact, $\B_{\aleph_0}\setminus\B_2$ is infinite -- see Proposition~\ref{counterex} below.

We begin with a useful definition. Let $x\sim(a_1,a_2,\dots)_q$; we say that $a_m$ is {\em forced} if there is no expansion of $x$ in base~$q$ of the form $x\sim(a_1,\dots,a_{m-1},b_m,\dots)_q$ with $b_m\neq a_m$.

\begin{figure}[t]
    \centering 
\centerline{
\begin{picture}(340,260)
\thicklines \put(-5,250){$x$}
\put(5,250){\line(1,0){5}}
\put(10,250){\line(1,0){300}} \put(10,200){\line(1,0){300}}
\put(60,150){\line(1,0){250}} \put(110,100){\line(1,0){200}}
\put(160,50){\line(1,0){150}}
\put(10,250){\line(0,-1){50}} \put(60,200){\line(0,-1){50}}
\put(110,150){\line(0,-1){50}} \put(160,100){\line(0,-1){50}}
\put(210,50){\line(0,-1){20}}
\put(315,248){$\ldots$} \put(315,198){$\ldots$}
\put(315,148){$\ldots$} \put(315,98){$\ldots$}
\put(315,48){$\ldots$} \put(209,15){$\vdots$}
\end{picture}}
\caption{A branching for countably many expansions}
    \label{fig:branching2}
\end{figure}
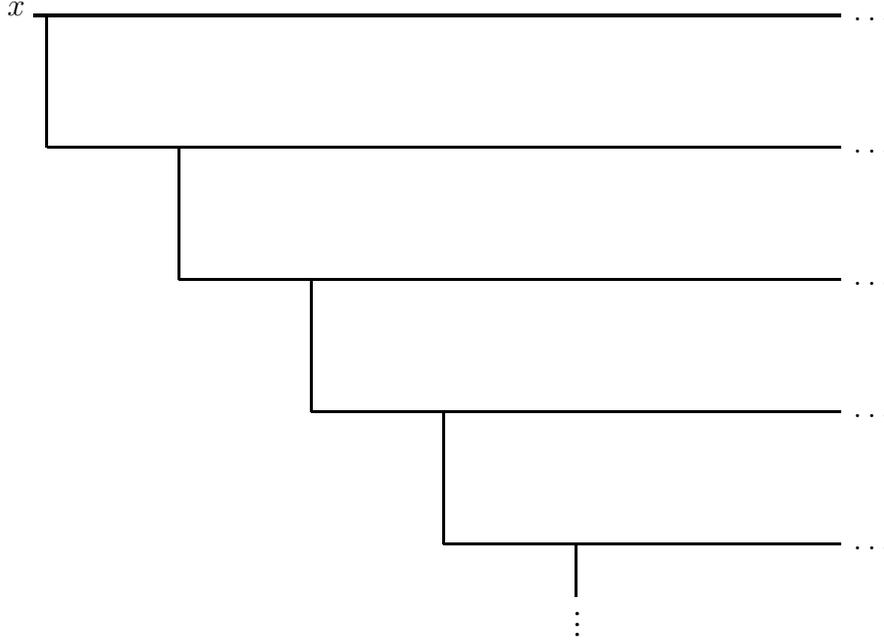

\begin{lemma}\label{lemaleph0}
Let $q>G$ and $x\sim(a_1,\dots,a_m,(01)^\infty)_q
\sim(b_1,\dots,b_k,a_1,\dots,a_m,(01)^\infty)_q$, where $a_1\neq b_1$, and assume that $a_2,\dots,a_m$ are forced in the first expansion and $b_2,\dots,b_k$ are forced in the second expansion. Then $q\in\B_{\aleph_0}$.
\end{lemma}
\begin{proof}Since all the symbols in the first expansion except $a_1$, are forced, the set of expansions for $x$ in base~$q$ is as follows:
\begin{align*}
&a_1,\dots,a_m,(01)^\infty, \\
&b_1,\dots,b_k,a_1,\dots,a_m,(01)^\infty, \\
&b_1,\dots,b_k,b_1,\dots,b_k,a_1,\dots,a_m,(01)^\infty, \\
&\vdots
\end{align*}
i.e., clearly infinite countable. The ``ladder'' branching pattern for $x$ is depicted in Fig.~\ref{fig:branching2}.
\end{proof}

\begin{prop}\label{counterex}
The set $\B_{\aleph_0}\cap(q_2,q_f)$ is infinite countable.
\end{prop}
\begin{proof}Define $q^{(n)}$ as the unique positive solution of
$$
(10000(10)^\infty)_{q^{(n)}}\sim(\boxed{0}\,11(01)^{n-1}
\boxed{1}\,0000(10)^\infty)_{q^{(n)}}
$$
(never mind the boxes for the moment) and put $\la_n=1/q^{(n)}$. A direct computation shows that
\[
\la_n^{2n+1}=\frac{1-\la_n-2\la_n^2+\la_n^3+\la_n^5}
{1-\la_n-2\la_n^2+\la_n^5},
\]
whence $\la_n\nearrow 1/q_2$ (as $1/q_2$ is a root of
$1-x-2x^2+x^3+x^5$), and consequently, $q^{(n)}\searrow q_2$.

By Lemma~\ref{lemaleph0}, if we show is that each symbol between the boxed 0 and the boxed 1 is forced, then $q^{(n)}\in\B_{\aleph_0}$. Let us prove it.

Notice that if $x\sim(a_1,a_2,\dots)_q$, then $a_1=0$ is forced if and only if $\sum_1^\infty a_kq^{-k}<1/q$; similarly, $a_1=1$ is forced if $\sum_1^\infty a_kq^{-k}>1/q(q-1)$. We need the following

\begin{lemma}\label{lem:forced}
\begin{enumerate}
\item If $q>G,m\ge0$ and $x\sim (1(01)^m1*)_q$, then the first 1 is forced (where $*$ stands for an arbitrary tail).
\item If $q>q_2, m\ge1$ and $x\sim((01)^m10000(10)^\infty)_q$, then
the first 0 is forced.
\end{enumerate}
\end{lemma}
\begin{proof}(1) By the above remark, we need to show that
\begin{equation*}
\frac1q+\frac1{q^3}+\dots+\frac1{q^{2m+1}}+\frac1{q^{2m+2}}>
\frac1{q(q-1)},
\end{equation*}
which is equivalent to (with $\la=1/q<1/G$)
\[
\frac{1-\la^{2m+2}}{1-\la^2}+\la^{2m+1}>\frac{\la}{1-\la}
\]
or $1-\la-\la^2>\la^{2m+1}-\la^{2m+2}-\la^{2m+3}$, which is true, in view of $1-\la-\la^2>0$ and $\la^{2m+1}<1$.

\medskip\noindent (2) Putting $\la=1/q$, we need to show that
\[
\la^2+\la^4+\dots+\la^{2m}+\la^{2m+1}+ \frac{\la^{2m+6}}{1-\la^2}<\la.
\]
This is equivalent to
\[
\la^{2m}<\frac{1-\la-\la^2}{1-\la-\la^2+\la^5}.
\]
The LHS in this inequality is a decreasing function of $m$, and for $m=1$ we have that it holds for $\la<0.59$, whence $q>q_2$ suffices.
\end{proof}

The proof of Proposition~\ref{counterex} now follows from the
definition of the sequence $(q^{(n)})_{n\ge1}$ and from
Lemma~\ref{lem:forced}.
\end{proof}

\begin{rmk}\label{rmk2}
The set $\B_{\aleph_0}\cap(G,q_2)$ is nonempty either: take $q_\om$ to be the appropriate root of $x^5=x^4+x^3+x-1$, with the numerical value $\approx1.68042$. Then
\[
x\sim(100(10)^\infty)_{q_\om}
\sim(\boxed{0}\,111\boxed{1}\,00(10)^\infty)_{q_\om},
\]
and similarly to the above, one can easily show that the three 1s between the boxed symbols are forced. Hence, by
Lemma~\ref{lemaleph0}, $q_\om\in\B_{\aleph_0}$. The question whether $\inf\B_{\aleph_0}=G$, remains open.
\end{rmk}

\begin{rmk} The condition of $q$ being transcendental in
Theorem~\ref{trans} is probably not necessary even for $q>q_f$. It would be interesting to construct an example of a family of
algebraic $q\in(q_f,\kl)$ for which the dichotomy in question holds.
\end{rmk}

\section{Middle order: $q$ just above $\kl$}

This case looks rather difficult for a hands-on approach, because, as we know, the set $\U_q$ for $q>\kl$ contains lots of transcendental numbers $x$, for which the tails of expansions in base~$q$ for $x$ and $x+1$ are completely different. However, a very simple argument allows us to link $\B_2$ to the well-developed theory of unique expansions for $x=1$.

Following \cite{KL}, we introduce
\[
\U:=\{q\in(1,2) : x=1\,\, \mbox{has a unique expansion in base
$q$}\}.
\]
Recall that in \cite{KL} it was shown that $\min\U=\kl$.

\begin{lemma}We have $\U\subset\B_2$. Consequently, the set
$\B_2\cap(\kl,\kl+\de)$ has the cardinality of the continuum for any $\de>0$.
\end{lemma}
\begin{proof}Since $x=0$ has a unique expansion in any base $q$, the first claim is a straightforward corollary of Lemma~\ref{crit}.

The second claim follows from the fact that $\U\cap(\kl,\kl+\de)$ has the cardinality of the continuum for any $\de>0$, which in turn is a consequence of the fact that the closure of $\U$ is a Cantor set -- see \cite[Theorem~1.1]{KL3}.
\end{proof}

\section{Top order: $q$ close to 2}\label{sec:top}

\subsection{$\bm{m=2}$.}
We are going to need the notion of thickness of a Cantor set. Our exposition will be adapted to our set-up; for a general case see, e.g., \cite{Ast}.

A Cantor set $C\subset\BR$ is usually constructed as follows: first we take a closed interval $I$ and remove a finite number of {\em gaps}, i.e., open subintervals of $I$. As a result we obtain a finite union of closed intervals; then we continue the process for each of these intervals {\em ad infinitum}. Consider the $n$th level, $\mathcal L_n$; we have a set of newly created gaps and a set of {\em bridges}, i.e., closed intervals connecting gaps. Each gap $\mathcal G$ at this level has two adjacent bridges, $\mathcal P$ and $\mathcal P'$.

The {\em thickness} of $C$ is defined as follows:
\[
\tau(C)=\inf_n\min_{\mathcal G\in\mathcal
L_n}\min\left\{\frac{|\mathcal P|}{|\mathcal G|}, \frac{|\mathcal
P'|}{|\mathcal G|}\right\},
\]
where $|I|$ denotes the length of an interval~$I$. For example, if $C$ is the standard middle-thirds Cantor set, then $\tau(C)=1$, because each gap is surrounded by two bridges of the same length.

The reason why we need this notion is the theorem due to Newhouse \cite{New} asserting that if $C_1$ and $C_2$ are Cantor sets, $I_1=\mbox{conv}(C_1), I_2=\mbox{conv}(C_2)$, and
$\tau(C_1)\tau(C_2)>1$ (where conv stands for convex hull), then $C_1+C_2=I_1+I_2$, provided the length of $I_1$ is greater than the length of the maximal gap in $C_2$ and vice versa. In particular, if $\tau(C)>1$, then $C+C=I+I$.

Notice that $\U_q$ is symmetric about the centre of $I_q$ --- because whenever $x\sim(a_1,a_2,\dots)_q$, one has $\frac1{q-1}-x\sim(1-a_1,1-a_2,\dots)_q$. Recall that Lemma~\ref{crit} yields the criterion $1\in\U_q-\U_q$ for $q\in\B_2$. Thus, we have $\U_q=1/(q-1)-\U_q$, whence $\U_q-\U_q=\U_q+\U_q-1/(q-1)$. Hence our criterion can be rewritten as follows:

\begin{equation}\label{sum}
q\in\B_2\Longleftrightarrow\frac{q}{q-1}\in\U_q+\U_q.
\end{equation}

It has been shown in \cite{GS} that the Hausdorff dimension of
$\U_q$ tends to 1 as $q\nearrow2$. Thus, one might speculate that for $q$ large enough, the thickness of $\U_q$ is greater than 1, whence by the Newhouse theorem, $\U_q+\U_q=2I_q$, which implies the RHS of (\ref{sum}).

However, there are certain issues to be dealt with on this way. First of all, in \cite{KV} it has been shown that $\U_q$ is not necessarily a Cantor set for $q>\kl$. In fact, it may contain
isolated points and/or be non-closed. This issue however is not really that serious because $\U_q$ is known to differ from a Cantor set by a countable or empty set \cite{KV}, which is negligible in our set-up.

A more serious issue is the fact that even if the Hausdorff
dimension of a Cantor set is close to 1, its thickness can be very small. For example, if one splits one gap by adding a very small bridge, the thickness of a resulting Cantor set will become very small as well! In other words, $\tau$ is not at all an increasing function with respect to inclusion.

Nonetheless, the following result holds:

\begin{lemma}\label{thick} Let $T$ denote the real root of $x^3=x^2+x+1$,
$T\approx1.83929$. Then
\begin{equation}\label{thicksum}
\U_q+\U_q=\left[0,\frac2{q-1}\right],\quad q\ge T.
\end{equation}
\end{lemma}
\begin{proof} Let $\Si_q$ denote the set of all sequences which provide unique expansions in base~$q$. It has been proved in
\cite{GS} that $\Si_q\subseteq\Si_{q'}$ if $q<q'$; hence
$\Si_T\subseteq\Si_q$. Note that by \cite[Lemma~4]{GS}, $\Si_T$ can be described as follows: it is the set of all 0-1 sequences which do not contain words $0111$ and $1000$ and also do not end with $(110)^\infty$ or $(001)^\infty$. Let $\widetilde\Si_T\supset\Si_T$ denote the set of 0-1 sequences which do not contain words $0111$ and $1000$. Note that by the cited lemma, $\widetilde\Si_T\subset\Si_q$ whenever $q>T$.

Denote by $\pi_q$ the projection map from $\{0,1\}^\BN$ onto $I_q$ defined by the formula
\[
\pi_q(a_1,a_2,\dots)=\sum_{n=1}^\infty a_nq^{-n},
\]
and put $\V_q=\pi_q(\widetilde\Si_T)$. Since $\widetilde\Si_T$ is a perfect set in the topology of coordinate-wise convergence, and since $\pi_q^{-1}|_{\U_q}$ is a continuous bijection, $\pi_q:\widetilde\Si_T\to\V_q$ is a homeomorphism, whence $\V_q$ is a Cantor set which is a subset of $\U_q$ for $q>T$. If $q=T$, then $\pi_q(\Si_T)=\pi_q(\widetilde\Si_T)$, hence the same conclusion about $\V_T$.

In view of Newhouse's theorem, to establish (\ref{thicksum}), it suffices to show that $\tau(\V_q)>1$, because
$\mbox{conv}(\V_q)=\mbox{conv}(\U_q)=I_q$. To prove this, we need to look at the process of creation of gaps in $I_q$. Note that any gap is the result of the words 000 and 111 in the symbolic space being forbidden. The first gap thus arises between $\pi_q([0110])=\bigl[\la^2+\la^3,\la^2+\la^3+
\frac{\la^5}{1-\la}\bigr]$ and $\pi_q([1001])=\bigl[\la+\la^4,\la+\la^4+\frac{\la^5}{1-\la}\bigr]$.
(Here, as above, $\la=q^{-1}\in(1/2,1/T)$ and $[i_1\dots i_r]$ denotes the corresponding cylinder in $\{0,1\}^\BN$.) The length of the gap is $\la+\la^4-\bigl(\la^2+\la^3+\frac{\la^5}{1-\la}\bigr)$, which is significantly less than the length of either of its adjacent bridges.

Furthermore, it is easy to see that any new gap on level~$n\ge5$ always lies between $\pi_q([a0110])$ and $\pi_q([a1001])$, where $a$ is an arbitrary 0-1 word of the length $n-4$ which contains neither 0111 nor 1000. The length of the gap is thus independent of $a$ and equals
$\la^{n-3}+\la^n-\la^{n-2}-\la^{n-1}-\frac{\la^{n+1}}{1-\la}$.

As for the bridges, to the right of this gap we have at least the union of the images of the cylinders $[a1001], [a1010]$ and $[a1011]$, which yields the length $\la^{n-3}+\frac{\la^{n-1}}{1-\la}-\la^{n-3}-\la^n=
\frac{\la^{n-1}}{1-\la} -\la^n$.

Hence
\begin{align*}
\frac{|\mbox{gap}|}{|\mbox{bridge}_1|}&\le\frac{1-\la-\la^2+\la^3-
\frac{\la^4}{1-\la}}{\frac{\la^2}{1-\la}-\la^3}\\
&=\frac{1-2\la+2\la^3-2\la^4}{\la^2-\la^3+\la^4}.
\end{align*}
This fraction is indeed less than 1, since this is equivalent to the inequality
\begin{equation}\label{ineq6}
3\la^4-3\la^3+\la^2+2\la-1>0,
\end{equation}
which holds for $\la>0.48$.

The bridge on the left of the gap is
$[\pi_q(a010^\infty),\pi_q(a01101^\infty)]$, and its length is
$|\mbox{bridge}_2|=
\la^{n-2}+\la^{n-1}+\frac{\la^{n+1}}{1-\la}-\la^{n-2}=
\la^{n-1}+\frac{\la^{n+1}}{1-\la}=\frac{\la^{n-1}}{1-\la}-\la^n=
|\mbox{bridge}_1|$, whence $|\mbox{bridge}_2|<|\mbox{gap}|$, and we are done.
\end{proof}

As an immediate corollary of Lemma~\ref{thick} and (\ref{sum}), we obtain

\begin{thm}\label{top}
For any $q\in [T,2)$ there exists $x\in I_q$ which has exactly two expansions in base~$q$.
\end{thm}

\begin{rmk} The constant $T$ in the previous theorem is clearly not sharp -- inequality~(\ref{ineq6}), which is the core of our proof, is essentially the argument for which we need a constant close to $T$. Considering $\U_q$ directly (instead of $\V_q$) should help decrease the lower bound in the theorem (although probably not by much).
\end{rmk}

\subsection{$\bm{m\ge3}$.}
\begin{thm}\label{close2}
For each $m\in\BN$ there exists $\ga_m>0$ such that
\[
(2-\ga_m,2)\subset \B_j,\quad 2\le j\le m.
\]
Furthermore, for any fixed $m\in\BN$,
\begin{equation}\label{eq:hdim}
\lim_{q\nearrow2}\dim_H \U_q^{(m)}=1,
\end{equation}
where, as above, $\U_q^{(m)}$ denotes the set of $x\in I_q$ which have precisely $m$ expansions in base~$q$.
\end{thm}
\begin{proof}Note first that if $q\in\B_m$ and
$1\in\U_q^{(m)}-\U_q$, then $q\in\B_{m+1}$. Indeed, analogously to the proof of Lemma~\ref{crit}, if $y\in\U_q$ and $y+1\in\U_q^{(m)}$, then $(y+1)/q$ lies in the interval $J_q$, and the shift of its expansion beginning with 1, belongs to $\U_q$, and the shift of its expansion beginning with 0, has $m$ expansions in base~$q$.

Similarly to Lemma~\ref{thick}, we want to show that for a fixed $m\ge2$,
\[
\U_q^{(m)}-\U_q=\left[-\frac1{q-1},\frac1{q-1}\right]
\]
if $q$ is sufficiently close to 2. We need the following result which is an immediate corollary of \cite[Theorem~1]{HKY}:

\medskip\noindent {\bf Proposition}. {\em For each $E>0$ there exists $\De>0$ such that for any two Cantor sets $C_1,C_2\subset\BR$ such that $\mathrm{conv}(C_1)=\mathrm{conv}(C_2)$ and $\tau(C_1)>\De,
\tau(C_2)>\De$, their intersection $C_1\cap C_2$ contains a Cantor set $C$ with $\tau(C)>E$.}

\medskip\noindent Let $T_k$ be the appropriate root of
$x^k=x^{k-1}+x^{k-2}+\dots+x+1$. Then $T_k\nearrow 2$ as
$k\to+\infty$, and it follows from \cite[Lemma~4]{GS} that
$\Si_{T_k}$ is a Cantor set of 0-1 sequences which do not contain $10^k$ nor $01^k$ and do not end with $(1^{k-1}0)^\infty$ or $(0^{k-1}1)^\infty$. Similarly to the proof of Lemma~\ref{thick}, we introduce the sets $\widetilde\Si_{T_k}$ and define
$\V_q^{(k)}=\pi_q(\widetilde\Si_{T_k})$ for $q>T_k$. For the same reason as above, $\V_q^{(k)}$ is always a Cantor set for $q\ge T_k$.

Using the same arguments as in the aforementioned proof, one can show that for any $M>1$ there exists $k\in\BN$ such that $\tau\bigl(\V_q^{(k)}\bigr)>M$ for all $q>T_k$. More precisely, any gap which is created on the $n$th level is of the form $[\pi_q(a01^{k-1}0)+\la^{n+1}/(1-\la),\pi_q(10^{k-1}1)]$, while the bridge on the right of this gap is at least $[\pi_q(10^{k-1}1),\pi_q(101^{k-1})+\la^{n+1}/(1-\la)]$; a simple computation yields
\begin{equation}\label{eq:tauinf}
\frac{|\mbox{bridge}|}{|\mbox{gap}|}\ge\frac{\la^2+\la^k-\la^{k+1}} {2\la^k-\la^{k+1}+1-2\la}\sim\frac{\la^2}{1-2\la}\to+\infty,\quad k\to+\infty,
\end{equation}
since $\la\le T_k^{-1}\to1/2$ as $k\to+\infty$. The same argument works for the bridge on the left of the gap.

We know that $\frac1q(\U_q^{(m)}\cap(\U_q-1))$ is a subset of
$\U_q^{(m+1)}\cap J_q$; we can also extend it to $I_q\setminus J_q$ by adding any number 0s or any number of 1s as a prefix to the expansion of any $x\in\frac1q(\U_q^{(m)}\cap(\U_q-1))$. Thus, $\mbox{conv}\bigl(\U_q^{(m+1)}\bigr)=I_q$, provided this set is nonempty.

Let us show via an inductive method that $\U_q^{(m+1)}$ is nonempty for $m\ge2$. Consider $\U_q^{(3)}$; by the above, there exists $k_3$ such that for $q>T_{k_3}$, the intersection
$\U_q\cap(\U_q-1)$ contains a Cantor set of thickness greater
than~1. Extending it to the whole of $I_q$, we obtain a Cantor set of thickness greater than 1 whose support is $I_q$. This set is contained in $\U_q^{(2)}$, whence
$\U_q^{(2)}-\U_q=\left[-\frac1{q-1},\frac1{q-1}\right]$, yielding that $\U_q^{(3)}\neq\emptyset$ for $q>T_{k_3}$.

Finally, by increasing $q$, we make sure $\U_q^{(3)}$ contains a Cantor set of thickness greater than 1, which implies
$\U_q^{(4)}\neq\emptyset$, etc. Thus, for any $m\ge3$ there exists $k_m$ such that $\U_q^{(m)}\neq\emptyset$ if $q>T_{k_m}$. Putting $\ga_m=2-T_{k_m}$ completes the proof of the first claim of the theorem.

To prove (\ref{eq:hdim}), notice that from (\ref{eq:tauinf}) it follows that $\tau(\V_q^{(k)})\to+\infty$ as $k\to+\infty$. Since for any Cantor set~$C$,
\[
\dim_H(C)\ge \frac{\log2}{\log(2+1/\tau(C))}
\]
(see \cite[p.~77]{PT}), we have $\dim_H(\V_q^{(k)})\to1$ as $k\to+\infty$, whence $\dim_H(\U_q^{(m)})\to1$ as $q\to2$, in view of $T_k\nearrow2$ as $k\to+\infty$.
\end{proof}

\begin{rmk} From the proof it is clear that the constructed sequence $\ga_m\to0$ as $m\to+\infty$. It would be interesting to obtain some bounds for $\ga_m$; this could be possible, since we roughly know how $\De$ depends on $E$ in the proposition quoted in the proof. Namely, from \cite[Theorem~1]{HKY} and the remark in p.~888
of the same paper, it follows that for {\em large} $E$ we have
$\De\sim\sqrt E$.
\end{rmk}

Finally, in view of $\ga_m\to0$, one may ask whether actually
$\bigcap\limits_{m\in\BN\cup\aleph_0}\B_m\neq\emptyset$. It turns out that the answer to this question is affirmative.

\begin{prop}\label{trib}
For $q=T$ and any $m\in\BN\cup\aleph_0$ there exists $x_m\in I_q$ which has $m$ expansions in base~$q$.
\end{prop}
\begin{proof} Let first $m\in\BN$. We claim that
\[
x_m\sim(1(000)^m(10)^\infty)_q\in\U_q^{(m+1)}.
\]
Note first that if some $x\sim(10001\dots)_q$ has an expansion
$(0,b_2,b_3,b_4,\dots)$ in base~$q$, then $b_2=b_3=b_4=1$ -- because $(01101^\infty)_q < (100010^\infty)_q$, a straightforward check.

Therefore, if $x_m\sim(0,b_2,b_3,\dots)_q$ for $m=1$, then
$b_2=b_3=b_4=1$, and since $1=1/q+1/q^2+1/q^3$, we have
$(b_5,b_6,\dots)_q=((10)^\infty)_q$, the latter being a unique
expansion. Hence $x_1$ has only two expansions in base~$q$.

For $m\ge2$, we still have $b_2=b_3=1$, but $b_4$ can be equal to 0. This, however, prompts $b_5=b_6=1$, and we can continue with $b_{3i-1}=b_{3i}=1,b_{3i+1}=0$ until $b_{3j-1}=b_{3j}=b_{3j+1}=1$ for some $j\le m$, since $(10^\infty)_q\sim((011)^j10^\infty)_q$,
whence $(1(000)^j(10)^\infty)_q>((011)^j01^\infty)_q$ for any
$j\ge1$.

Thus, any expansion of $x_m$ in base~$q$ is of the form
\[
x_m\sim((011)^j1(000)^{m-j}(10)^\infty)_q,\quad 0\le j\le m,
\]
i.e., $x_m\in\U_q^{(m+1)}$.

For $m=\aleph_0$, we have
$x_\infty=\lim_{m\to\infty}x_m\sim(10^\infty)_q$. Notice that
$x_\infty\sim(011\,\,10^\infty)_q$ with the first two 1s clearly forced so we can apply Lemma~\ref{lemaleph0} to conclude that
$x_\infty\in\U_q^{(\aleph_0)}$, whence $q=T\in\B_{\aleph_0}$ as well.
\end{proof}

\begin{rmk} The choice of the tail $(10)^\infty$ in the proof
is unimportant; we could take any other tail, as long as it is a unique expansion which begins with 1. Thus, for $q=T$,
\[
\dim_H\U_q^{(m)}=\dim_H\U_q,\quad m\in\BN.
\]
This seems to be a very special case, because typically one might expect a drop in dimension with $m$. Note that in \cite{GS} it has been shown that $\dim_H\U_T=\log G/\log T\approx0.78968$.
\end{rmk}

\section{Summary and open questions}\label{final}

Summing up, here is the list of basic properties of the set $\B_2$:
\begin{itemize}
\item The set $\B_2\cap(G,\kl)$ is infinite
countable and contains only algebraic numbers (the ``lower
order''\footnote{Our terminology is borrowed from cricket.}). The latter claim is valid for $\B_m$ with $m\ge3$, although it is not clear whether $\B_m\cap(G,\kl)$ is nonempty.
\item $\B_2\cap(\kl,\kl+\de)$ has the cardinality of the continuum for any $\de>0$ (the ``middle order'').
\item $[T,2)\subset\B_2$ (the ``top order''), with a similar claim about $\B_m$ with $m\ge3$.
\end{itemize}

Here are a few open questions:

\begin{itemize}
\item Is $\B_2$ closed?
\item Is $\B_2\cap(G,\kl)$ a discrete set?
\item Is it true that $\dim_H(\B_2\cap(\kl,\kl+\de))>0$ for any
$\de>0$?
\item Is it true that $\dim_H(\B_2\cap(\kl,\kl+\de))<1$ for
some $\de>0$?
\item What is the value of $\inf\B_m$ for $m\ge3$?
\item What is the smallest value $q_0$ such that
$\U_q+\U_q=2I_q$ for $q\ge q_0$?
\item Is $\inf \B_{\aleph_0}=G$?
\item Does $\B_{\aleph_0}$ contain an interval as well?
\end{itemize}

\medskip\noindent\textbf{Acknowledgement.} The author is grateful to Boris Solomyak for stimulating questions and to Martijn de Vries and especially the anonymous referee for indicating some errors in the preliminary version of this paper.

\begin{table}[ht]
\centerline{
\begin{tabular}{c|c|c}
       $q$ & equation & numerical value \\
\hline $G$ & $x^2=x+1$ & 1.61803 \\
\hline $q_\om$ & $x^5=x^4+x^3+x-1$ & 1.68042 \\
\hline $q_2$ & $x^4=2x^2+x+1$ & 1.71064 \\
\hline $q_f$ & $x^3=2x^2-x+1$ & 1.75488 \\
\hline $\kl$ & $\sum_{1}^{\infty}\mathfrak{m}_{n}x^{-n+1}=1$ & 1.78723 \\
\hline $T$ & $x^3=x^2+x+1$ & 1.83929 
\end{tabular}
}

\vskip0.5truecm

\caption{The table of constants used in the text.} \label{table}
\end{table}

\end{document}